\documentclass[a4paper,12pt]{article} 
\usepackage{amssymb}
\usepackage{amscd}
\usepackage{amsmath} 
\usepackage[dvipdfmx]{graphicx}
\usepackage{latexsym} 
\usepackage{enumerate}
\usepackage[usenames]{color}
\usepackage{lscape}

\usepackage{theorem}
\theoremstyle{plain}
\theorembodyfont{\itshape}
\newtheorem{theorem}{Theorem}[section]
\newtheorem{proposition}[theorem]{Proposition}
\newtheorem{lemma}[theorem]{Lemma}
\newtheorem{corollary}[theorem]{Corollary}

\theorembodyfont{\rmfamily}

\newtheorem{remark}[theorem]{Remark}

\newtheorem{example}[theorem]{Example}
\newtheorem{assumption}[theorem]{Assumption}

\def\re{\mathrm{e}}

\def\Mm{\mathrm{M{\scriptstyle ax}}}

\def\dist{\mathrm{dist}}

\def\rRe{\mathrm{Re}} 
 
\def\bC{\mathbb{C}}

\def\bN{\mathbb{N}}

\def\bQ{\mathbb{Q}}
\def\bR{\mathbb{R}}
\def\bZ{\mathbb{Z}}

\def\cE{\mathcal{E}}

\def\ve{\varepsilon}
\def\vG{\varGamma}
\def\vD{\varDelta}

\def\ba{\mbox{\boldmath $a$}}

\def\bb{\mbox{\boldmath $b$}}

\def\bk{\mbox{\boldmath $k$}}

\def\bp{\mbox{\boldmath $p$}}

\def\bal{\mbox{\boldmath $\alpha$}}
 
\def\blambda{\mbox{\boldmath $\lambda$}} 
\def\bmu{\mbox{\boldmath $\mu$}} 
\def\bsigma{\mbox{\boldmath $\sigma$}} 
\def\btau{\mbox{\boldmath $\tau$}} 
 
\def\0{\mbox{\boldmath $0$}}
\def\1{\mbox{\boldmath $1$}}
\def\2{\mbox{\boldmath $2$}}
\def\3{\mbox{\boldmath $3$}}
\def\cfL{\mbox{\Large \boldmath $\mathrm{K}$}}

\def\Mm{\mathrm{M{\scriptstyle ax}}}

\def\hgf{{}_2f_1}
\def\hgF{{}_2F_1}

\title{\bf Discrete Laplace Method and Truncation \\ 
Error of Gauss Continued Fraction}   
\author{
Katsunori Iwasaki\thanks{Department of Mathematics, 
Hokkaido University, Kita 10, Nishi 8, Kita-ku, Sapporo 060-0810 Japan. 
{\tt iwasaki@math.sci.hokudai.ac.jp}}} 
\date{April 6, 2019}  
\begin{document}
\maketitle
\begin{abstract} 
The leading asymptotics of the truncation 
error for Gauss's continued fraction is determined exactly.    
Not only for this purpose but also for wider applicability elsewhere 
the discrete analogue of Laplace's method for hypergeometric series 
containing a large parameter, which was developed in a previous 
paper,  is generalized  in two directions.         
\end{abstract} 
\section{Introduction} \label{sec:intro}
In 1813 Gauss introduced a general continued fraction   
\begin{equation} \label{eqn:cfG}
\overset{\infty}{\underset{n=0}{\cfL}} \,\, \frac{R(n)}{1} 
= \frac{R(0)}{1} \,\, \mathop{}_{+} \,\, \frac{R(1)}{1} 
\,\, \mathop{}_{+} \,\, \frac{R(2)}{1} \, \, \mathop{}_{+} \,\, 
\mathop{}_{\textstyle \cdots,}   
\end{equation}
known today as Gauss's continued fraction (GCF for short), where 
$R(0) := 1$ and 
\begin{align*}
R(2 m+1) &:= -\frac{(m+b)(m+c-a)z}{(2m+c)(2m+c+1)},  \\[1mm]
R(2 m+2) &:= -\frac{(m+a+1)(m+c-b+1)z}{(2m+c+1)(2m+c+2)},    
\qquad m \in \bZ_{\ge 0},   
\end{align*}
where $a$, $b$, $c$ and $z$ are complex parameters, with $z$ being 
referred to as the independent variable. 
For non-vanishing of the numerators and denominators of 
$R(n)$, $n \in \bZ_{\ge 1}$, we assume that 
\begin{equation} \label{eqn:abc}
a, \,\, c-b \, \not\in \, \bZ_{\le -1}; \qquad 
b, \,\, c, \,\, c-a \, \not\in \, \bZ_{\le 0}. 
\end{equation}
It is well known that for $z \in \bC_z \setminus [1, \, \infty )$ the 
continued fraction \eqref{eqn:cfG} converges to the ratio 
$$
\frac{\hgF(a+1, b; c+1; z)}{ \hgF(a, b; c; z)},  
$$
where $\hgF(a, b; c; z)$ represents Gauss's hypergeometric 
series as well as its analytic continuation to the cut plane 
$\bC \setminus [1, \, \infty)$; 
see e.g. Jones and Thron \cite[Theorem 6.1]{JT}.    
\par
Let $\ba := (a, b; c)$, $\bk := (1, 0; 1)$ and $\bp := \bk + \sigma(\bk) = (1,1;2)$, 
where $\sigma : (a, b; c) \mapsto (b, a; c)$ exchanges the upper  
parameters $a$ and $b$.   
Notice that $\hgF(\ba; z)$ is invariant under the involution $\sigma$.  
Continued fraction \eqref{eqn:cfG} is associated with three-term 
contiguous relations
\begin{subequations} \label{eqn:Contig}
\begin{align} 
\hgF(\ba; z) &= \hgF(\ba+\bk; z) - \frac{b(c-a)z}{c(c+1)} \, \hgF(\ba+\bp; z),   
\label{eqn:Contig1} \\  
\hgF(\ba + \bk; z) &= \hgF(\ba+\bp; z) - 
 \frac{(a+1)(c-b+1)z}{(c+1)(c+2)} \, \hgF(\ba+\bp+\bk; z),  
\label{eqn:Contig2}  
\end{align}
\end{subequations}
where \eqref{eqn:Contig1} can be found in Andrews et al. 
\cite[formula (2.5.11)]{AAR},   
while \eqref{eqn:Contig2} is obtained from \eqref{eqn:Contig1} by applying 
$\sigma$ and replacing $\ba$ with $\ba + \bk$.  
For $m \in \bZ_{\ge 0}$ let
$$
F(2m) := \hgF(\ba+ m \bp; z), \qquad F(2m+1) := \hgF(\ba+ m \bp + \bk; z).  
$$
Taking shifts $\ba \mapsto \ba + m \bp$ in \eqref{eqn:Contig} 
induces a three-term recurrence relation 
\begin{equation} \label{eqn:Rr}
F(n) = F(n+1) + R(n+1) \, F(n+2), \qquad n \in \bZ_{\ge 0},   
\end{equation}
where $n$ is either $2m$ or $2m+1$.  
Continued fraction \eqref{eqn:cfG} then follows from \eqref{eqn:Rr} formally. 
\par
We are interested in the truncation error of Gauss's continued fraction,   
\begin{equation} \label{eqn:En} 
\cE_n(\ba; z) := \frac{\hgF(a+1, b; c+1; z)}{ \hgF(a, b; c; z)} - 
\overset{n}{\underset{j=0}{\cfL}} \,\, \frac{R(j)}{1}.    
\end{equation}
It is also interesting to consider the specialization of letting $a \to 0$ 
followed by the substitution $c \mapsto c-1$. 
The truncation error \eqref{eqn:En} then turns into   
$$
\cE_n^*(\bb; z) := \hgF(1, b; c; z) - 
\overset{n}{\underset{j=0}{\cfL}} \,\, \frac{R^*(j)}{1},  \qquad \bb := (b; c),  
$$
where $R^*(0) := 1$ and $R^*(n)$ with $n = 2m+1$ or $2m+2$ is given by 
\begin{align*}
R^*(2 m+1) &:= -\frac{(m+b)(m+c-1)z}{(2m+c-1)(2m+c)},  \\[1mm]
R^*(2 m+2) &:= -\frac{(m+1)(m+c-b)z}{(2m+c)(2m+c+1)},    
\qquad m \in \bZ_{\ge 0}.    
\end{align*}
In order for $R^*(n)$, $n \in \bZ_{\ge 1}$, not to be indefinite, 
we assume that 
\begin{equation} \label{eqn:bc}
b, \,\, c, \,\, c-b \not\in \bZ_{\le 0}. 
\end{equation}
\par
 J.~Borwein et al. \cite[Theorem 4]{BCP} gave the following estimate 
in a special case of Gauss's continued fraction: 
If $(\bb, z)$ satisfies $2 \le b$, $b+1 \le c \le 2 b$ and 
$-1 \le z < 0$, then 
$$
|\cE^*_n(\bb; z)|  \le 
\frac{\vG(m+1) \, (m+b) \, \vG(m+c-b) \, \vG(b) \vG(c)}{\vG(m+b) \, 
\vG(m+c) \, b \, \vG(c-b)} 
\left\{ \frac{2 b}{(c-2) \left(1-\frac{2}{z} \right) + (2 b-c)} \right\}^n, 
$$
where $m := \lfloor n/2 \rfloor$ is the largest integer not exceeding 
$n/2$. 
As another topic, based on Gauss's continued fraction and other means, 
Colman et al. \cite{CCD} developed an efficient algorithm for the validated 
high-precision computation of certain special $\hgF$ functions.   
\par
The purpose of this article is to determine the leading asymptotics 
of the truncation error $\cE_n(\ba; z)$ as $n \to \infty$ for general 
$\ba = (a, b; c) \in \bC^3$ and $z \in (-\infty, \, 1)$. 
Given two sequences $\{\alpha_n\}$ and $\{\beta_n\}$, we mean by 
$\alpha_n \sim \beta_n$ that their ratio behaves like 
$\alpha_n/\beta_n = 1 + O(n^{-\frac{1}{2}})$ as 
$n \to \infty$. 
Then our main result is stated in the following manner. 
\begin{theorem} \label{thm:Error} 
If $(\ba; z)$ satisfies condition \eqref{eqn:abc}, $z \in (-\infty, \, 1)$ 
and $\hgF(\ba; z) \neq 0$, then 
\begin{equation} \label{eqn:Error}
\begin{split}
\cE_n(\ba; z) &\sim \frac{2 \pi }{\hgF(\ba; z)^2} \cdot 
\frac{\vG(c) \vG(c+1)}{\vG(a+1) \vG(b) \vG(c-a) \vG(c-b+1)} \\[1mm]
&\phantom{==} \times 
\frac{z (1-z)^{c-a-b}}{(1+\sqrt{1-z})^{2(c+1)}}
\left\{ \frac{z}{(1+\sqrt{1-z})^2}\right\}^{n}.     
\end{split} 
\end{equation}
\end{theorem}
\par
The relation $\sim$ in \eqref{eqn:Error} is compatible with the 
specialization and we have the following. 
\begin{corollary} \label{cor:Error} 
If $(\bb; z)$ satisfies condition \eqref{eqn:bc} and 
$z \in (-\infty, \, 1)$, then 
\begin{equation} \label{eqn:Error2}
\cE_n^*(\bb; z) \sim  
\frac{2\pi  \vG(c)}{\vG(b) \vG(c-b)} \cdot 
\frac{z (1-z)^{c-b-1}}{(1+\sqrt{1-z})^{2 c}} 
\left\{ \frac{z}{(1+\sqrt{1-z})^2}\right\}^{n}.  
\end{equation}
\end{corollary}
\par
For every $z \in (-\infty, \, 1)$ the dilation constant 
$z (1+\sqrt{1-z})^{-2}$ in \eqref{eqn:Error} and \eqref{eqn:Error2} is 
smaller than $1$ in its absolute value, so that $\cE_n(\ba; z)$ and 
$\cE_n^*(\bb; z)$ decay exponentially as $n \to \infty$.    
\par
Besides its intrinsic interest, the error estimate of Gauss's continued 
fraction is instructive as a testing ground for our discrete analogue 
of Laplace's method for general hypergeometric series containing 
a large parameter.  
The latter content is expected to have many applications to 
hypergeometric series, especially to those of higher order.   
Indeed, an earlier version of it has already had an interesting 
application to ${}_3F_2(1)$ continued fractions in \cite{EI2}.  
\par
In general a continued fraction is associated with a 
three-term recurrence relation and the truncation error 
of the former can be controlled by the ratio of a recessive 
sotution to a dominant one of the latter. 
For a hypergeometric continued fraction the associated 
recurrence relation comes from a contiguous relation.
For an efficient treatment of recessive and dominant 
solutions the contiguous relation should be rescaled 
in an appropriate sense.  
This is the theme of  ``simultaneous contiguous relations''  
in \S\ref{sec:scr}.  
In accordance with this rescaling, the rescaled Gauss continued fraction 
(rGCF for short) is introduced and its relation with the original GCF is 
established in \S \ref{sec:r-gcf}. 
Then the recurrence relation associated with the rGCF 
is considered.  
An asymptotic representation of a recessive solution to 
it is given in \S \ref{sec:clm}.  
\par
To deal with dominant solutions, we turn our attention to the 
general theory of discrete Laplace method.    
In \S \ref{sec:dlm} two improvements of the earlier version in 
\cite{EI2} are made to facilitate its broader applicability. 
This generalization is illustrated by a couple of examples in 
\S \ref{sec:ex}, which are chosen in anticipation of a later 
application to the rGCF. 
The assumption imposed in \S \ref{sec:dlm} is not always 
fulfilled by a general hypergeometric series.   
To cope with this situation one has to cut the series 
into several pieces and manipulate them so that the desired 
assumption is recovered for each component.  
The recipe for this procedure is given in \S\ref{sec:pi-sc}. 
In \S\ref{sec:dom} we return to the situation of rGCF and 
derive asymptotic formulas for two dominant solutions to the 
associated recurrence.  
In \S \ref{sec:error}, after calculating the Casoratian of 
recessive and dominant solutions, we establish 
Theorem \ref{thm:Error} and Corollary \ref{cor:Error} by 
putting all the discussions together.        
\section{Simultaneous Contiguous Relations} \label{sec:scr}
Consider a rescaled version of Gauss's hypergeometric series 
\begin{equation} \label{eqn:hgf}
\hgf(\ba; z) 
:= \sum_{k=0}^{\infty} \frac{\vG(a+k) \vG(b+k)}{\vG(1+k) \vG(c+k)} \, z^k 
= \frac{\vG(a) \vG(b)}{\vG(c)} \, \hgF(\ba; z).    
\end{equation}
For generic values of the parameters $\ba = (a, b; c) \in \bC^3$ 
we also consider the rescaled version of Frobenius solutions to 
the Gauss hypergeometric equation, 
\begin{align} 
f_1^{(0)}(\ba; z) &:= \hgf(a, \, b; \, c; \, z), 
\tag{E1} \label{E1} \\
f_2^{(0)}(\ba; z) &:= z^{1-c} \, \hgf(a-c+1, \, b-c+1; \, 2-c; \, z), 
\tag{E17} \label{E17} \\
f_1^{(1)}(\ba; z) &:= \hgf(a, \, b; \, a+b-c+1; \, 1-z), 
\tag{E5} \label{E5} \\
f_2^{(1)}(\ba; z) &:= (1-z)^{c-a-b} \, \hgf(c-a, \, c-b; \, c-a-b+1; \, 1-z), 
\tag{E21} \label{E21} \\
f_1^{(\infty)}(\ba; z) &:= (-z)^{-a} \, \hgf(a, \, a-c+1; \, a-b+1; \, z^{-1}), 
\tag{E9} \label{E9} \\
f_2^{(\infty)}(\ba; z) &:= (-z)^{-b} \, \hgf(b-c+1, \, b; \, b-a+1; \, z^{-1}).  
\tag{E13} \label{E13}
\end{align}
where for example \eqref{E17} indicates that the original non-rescaled 
solution appears as formula (17) in 
Erd\'{e}lyi et al. \cite[Chap. II, \S2.8]{Erdelyi}. 
It will be more convenient to take a further rescaling   
\begin{subequations} \label{eqn:r-frob}
\begin{alignat}{2}
y_1^{(0)}(\ba; z) &:= f_1^{(0)}(\ba; z), \qquad & 
y_2^{(0)}(\ba; z) &:= f_2^{(0)}(\ba; z),  \label{eqn:r-frob0} \\[1mm]
y_1^{(1)}(\ba; z) &:= \chi(\ba) \, f_1^{(1)}(\ba; z), \qquad & 
y_2^{(1)}(\ba; z) &:=  \chi(\ba) \, f_2^{(1)}(\ba; z), \label{eqn:r-frob1} 
\\[2mm] 
y_1^{(\infty)}(\ba; z) &:= \frac{f_1^{(\infty)}(\ba; z)}{\sin \pi(c-b)}, 
\qquad & 
y_2^{(\infty)}(\ba; z) &:= \frac{f_2^{(\infty)}(\ba; z)}{\sin \pi(c-a)}, 
\label{eqn:r-frob-inf}
\end{alignat}
\end{subequations}
where the multiplicative factor $\chi(\ba)$ in \eqref{eqn:r-frob1} 
is given by 
\begin{equation} \label{eqn:m1(a)}
\chi(\ba) := \frac{ \pi \sin \pi c}{\sin \pi(c-a) \cdot \sin \pi(c-b)} \cdot 
\frac{1}{\vG(c-a) \vG(c-b)}.  
\end{equation}
\par
The connection formulas for the rescaled Frobenius solutions 
\eqref{eqn:r-frob} are given by  
\begin{align}
y_1^{(1)}(\ba; z) 
&= y_1^{(0)}(\ba; z) - y_2^{(0)}(\ba; z), \tag{E35} \label{E35}
\\[2mm]
y_2^{(1)}(\ba; z) 
&= \frac{\sin \pi a \cdot \sin \pi b }{\sin \pi(c-a) \cdot \sin \pi(c-b)} \, 
y_1^{(0)}(\ba; z) - y_2^{(0)}(\ba; z),  \tag{E43} \label{E43}
\\[2mm]
y_1^{(\infty)}(\ba; z)  
&= \frac{\sin \pi b}{\sin \pi c \cdot \sin \pi(c-b)} \, y_1^{(0)}(\ba; z) +  
 \frac{e^{i \pi c}}{\sin \pi c} \, y_2^{(0)}(\ba; z), \tag{E37} \label{E37} 
\\[2mm]
y_2^{(\infty)}(\ba; z) 
&= \frac{\sin \pi a}{\sin \pi c \cdot \sin \pi(c-a)} \, y_1^{(0)}(\ba; z) + 
\frac{e^{i \pi c}}{\sin \pi c} \, y_2^{(0)}(\ba; z),  \tag{E39} \label{E39} 
\end{align} 
where for example \eqref{E43} indicates that the original 
non-rescaled version can be found in formula (43) of Erd\'elyi et al. 
\cite[Chap. II, \S2.8]{Erdelyi}. 
It is remarkable that all of the rescaled connection coefficients are 
$\bZ^3$-{\sl periodic}, that is, invariant under the translation of 
$\ba$ by any integer vector.  
\par
For any nonzero integer vectors $\bk$, $\bp \in \bZ^3$ with $\bk \neq \bp$ 
there exist unique rational functions $u(\ba; z)$, $v(\ba; z) \in \bQ(\ba, z)$ 
such that $y_1^{(0)}(\ba; z) = \hgf(\ba; z)$ satisfies three-term relation 
\begin{equation} \label{eqn:contig}
y(\ba; z) = u(\ba; z) \, y(\ba+\bk; z) + v(\ba; z) \, y(\ba+\bp; z). 
\end{equation}
An equation of this sort is called a {\sl contiguous relation}. 
An argument in \cite[\S 2]{EI1} (which deals with ${}_3F_2(1)$ 
but remains valid for $\hgF$) shows that the other rescaled Frobenius 
solution $y_2^{(0)}(\ba; z)$ at the origin satisfies the same contiguous 
relation \eqref{eqn:contig}.  
It then follows from the connection formulas mentioned above, especially 
from the $\bZ^3$-periodicity of the connection coefficients, that 
contiguous relation \eqref{eqn:contig} is satisfied by all the six 
rescaled Frobenius solutions \eqref{eqn:r-frob}.  
We refer to this property as the {\sl simultaneousness} of contiguous 
relations. 
\section{Rescaled Gauss Continued Fraction} \label{sec:r-gcf}
The simultaneous contiguous relations corresponding to 
\eqref{eqn:Contig1} and \eqref{eqn:Contig2} are given by 
\begin{subequations} \label{eqn:contig2}
\begin{align}  
y(\ba; z) 
&= \frac{c}{a} \, y(\ba+\bk; z) + \frac{(a-c)z}{a} \, y(\ba+\bp; z),   
\label{eqn:contig21}  \\[1mm] 
y(\ba + \bk; z) 
&= \frac{c+1}{b} \, y(\ba+\bp; z) +\frac{(b-c-1)z}{b} \, 
y(\ba+\bp+\bk; z),   \label{eqn:contig22} 
\end{align}
\end{subequations} 
where $y(\ba; z)$ is any member of the six functions in \eqref{eqn:r-frob}  
and $\bk := (1, 0; 1)$, $\bp := \bk + \sigma(\bk) = (1,1;2)$ as in 
\S \ref{sec:intro}. 
For $m \in \bZ_{\ge 0}$ let $y(2 m) := y(\ba+ m \bp; z)$, 
$y(2m+1) := y(\ba+ m \bp + \bk; z)$ and  
\begin{alignat*}{2}
q(2m) &:= \frac{2 m+c}{m+a}, \qquad & r(2m+1) &:= -\frac{(m+c-a)z}{m+a}, 
\\[1mm]
q(2m+1) &:= \frac{2m+c+1}{m+b}, \qquad & r(2m+2) &:= -\frac{(m+c-b+1)z}{m+b}.  
\end{alignat*}
Taking shifts $\ba \mapsto \ba + m \bp$, $m \in \bZ_{\ge0}$ in 
\eqref{eqn:contig2} leads to a three-term recurrence relation 
\begin{equation} \label{eqn:rr}
y(n) = q(n) \, y(n+1) + r(n+1) \, y(n+2), \qquad n \in \bZ_{\ge 0}.  
\end{equation}
where $n$ is either $2 m$ or $2m+1$. 
If $y(\ba; z)$ is $y_i^{(*)}(\ba; z)$ in \eqref{eqn:r-frob} 
then $y(n)$ is denoted by $y_i^{(*)}(n)$. 
\begin{remark} \label{rem:Pfaff} 
Recall that there are two transformation formulas called Pfaff's  
transformations,  
$$
\hgF(a,b;c; z) 
= (1-z)^{-a} \hgF(a,c-b; c; z/(z-1)) 
= (1-z)^{-b} \hgF(c-a, b; c; z/(z-1)), 
$$
together with their composite called Euler's transformation (see 
e.g. \cite[Theorem 2.2.5]{AAR}). 
We can then speak of the rescaled version of these transformations  
for $y_i^{(*)}(\ba; z)$ and $y_i^{(*)}(n)$.    
\end{remark}
\par
Recurrence relation \eqref{eqn:rr} formally induces a rescaled version 
of Gauss's continued fraction 
\begin{equation} \label{eqn:cfg}
\overset{\infty}{\underset{n=0}{\cfL}} \,\, \frac{r(n)}{q(n)} 
= \frac{r(0)}{q(0)} \,\, \mathop{}_{+} \,\, \frac{r(1)}{q(1)} 
\,\, \mathop{}_{+} \,\, \frac{r(2)}{q(2)} \, \, \mathop{}_{+} \,\, 
\mathop{}_{\textstyle \cdots} \qquad \mbox{with} \quad 
r(0) := 1.     
\end{equation}
Continued fractions \eqref{eqn:cfG} and 
\eqref{eqn:cfg} are equivalent up to a constant multiple, 
more precisely, 
\begin{equation} \label{eqn:equiv}
\overset{n}{\underset{j=0}{\cfL}} \,\, \frac{R(j)}{1} 
=  
\frac{c}{a} \,\, \overset{n}{\underset{j=0}{\cfL}} \,\, \frac{r(n)}{q(n)}, 
\qquad n \in \bZ_{\ge 0}.  
\end{equation}
It will turn out that if $y_1^{(0)}(\ba; z) = \hgf(\ba; z)$ is chosen for $y(\ba; z)$,  
then the corresponding sequence $f(n) := y_1^{(0)}(n)$ is a 
{\sl recessive} solution to the recurrence equation \eqref{eqn:rr}. 
So Pincherle's theorem \cite[Theorem 5.7]{JT} implies that continued fraction 
\eqref{eqn:cfg} converges to the ratio 
$f(1)/f(0) = \hgf(\ba+\bk; z)/ \hgf(\ba; z)$. 
We are interested in the asymptotic behavior of the truncation error 
\begin{equation} \label{eqn:en}
\ve_n(\ba; z) := \frac{\hgf(\ba+\bk; z)}{ \hgf(\ba; z)} - 
\overset{n}{\underset{j=0}{\cfL}} \,\, \frac{r(j)}{q(j)} = 
\frac{c}{a} \, \cE_n(\ba; z).   
\end{equation}
where the second equality follows from definitions \eqref{eqn:En} and 
\eqref{eqn:hgf} and relation \eqref{eqn:equiv}. 
\par
If $g(n)$ is a {\sl dominant} solution to \eqref{eqn:rr} then the error estimate in 
\cite[\S 3.1, formula (29)]{EI2} reads 
\begin{equation} \label{eqn:error} 
\ve_n(\ba; z) = \frac{\omega(0) \cdot h(n)}{f(0)^2} 
\left\{ 1 + O\left(\frac{g(0) \cdot h(n)}{f(0)} \right)\right\} 
\qquad \mbox{as} \quad n \to \infty, 
\end{equation}
where $h(n) := f(n+2)/g(n+2)$ is the ratio of the recessive solution to the 
dominant one, while $\omega(n) := f(n) \cdot g(n+1) - f(n+1) \cdot g(n)$ is 
the Casoratian of $f(n)$ and $g(n)$. 
We remark that Landau's symbol in \eqref{eqn:error} is locally uniform 
with respect a parameter contained in it, so that even if $f(0)$, $g(0)$ 
and/or  $h(n)$ are individually singular at some value of the parameter, 
it remains valid as far as the expression $g(0) \cdot h(n)/f(0)$ 
is regular in total.  
\section{Recessive Solution} \label{sec:clm} 
Using the usual (continuous) Laplace method we shall find the asymptotic 
behavior of the sequence $y_1^{(0)}(n)$, which will serve as a 
recessive solution to the recurrence equation \eqref{eqn:rr}. 
Given two sequences $\{ \alpha_n \}$ and $\{ \beta_n \}$, we mean by 
$\alpha_n \approx \beta_n$ that  
$\alpha_n/\beta_n = 1 + O(n^{-1})$ as $n \to \infty$.  
\begin{proposition} \label{prop:clm} 
For any $z \in (-\infty, \, 1)$ there exists an asymptotic representation 
\begin{equation} \label{eqn:clm}
y_1^{(0)}(n) \approx \frac{2 \sqrt{\pi} \, (2 \sqrt{1-z})^{c-a-b-\frac{1}{2}}}{ 
n^{c-a-b+\frac{1}{2}} \, (1+\sqrt{1-z})^{n+c-1}}.   
\end{equation}
\end{proposition} 
{\it Proof}. 
For $\rRe \, a > 0$ and $\rRe(c-a) > 0$, Euler's integral representation reads  
$$
\hgf(\ba; z) = \frac{\vG(b)}{\vG(c-a)} \int_0^1 u(x) \, d x, 
\qquad u(x) := x^{a-1} (1-x)^{c-a-1} (1- z x)^{-b}.    
$$
Thus for $m \in \bZ_{\ge 0}$, $\rRe \, a > -m$ and $\rRe(c-a) > -m$, 
we have 
$$
\hgf(\ba + m \bp; z) = \frac{\vG(m+b)}{\vG(m+c-a)} 
\int_0^1 \Phi(x)^m \, u(x) \, d x, 
\qquad \Phi(x) := \frac{x(1-x)}{1- z x}.     
$$
The gamma factor behaves like $\vG(m+b)/\vG(m+c-a) \approx m^{a+b-c}$ 
due to Stirling's formula. 
\par
Define a function $\phi(x)$ by $\Phi(x) = \re^{- \phi(x)}$ and observe that  
$$
\phi'(x) = \frac{\phi_1(x)}{x(1-x)(1-z x)} 
\qquad \mbox{with} \qquad  \phi_1(x) := - z x^2 + 2 x -1. 
$$
The quadratic equation $\phi_1(x) = 0$ has a unique root 
$x_0 := (1 + \sqrt{1-z})^{-1}$ such that $0 < x_0 < 1$. 
Some calculations show that $\Phi(x_0) = 
(1+\sqrt{1-z})^{-2}$ and   
$$
\phi''(x_0) = \frac{2(1+\sqrt{1-z})^2}{\sqrt{1-z}} > 0,  
\qquad 
u(x_0) = \frac{\sqrt{1-z}^{\, c-a-b-1}}{(1+\sqrt{1 - z})^{c-2}}.    
$$
\par
A standard argument in Laplace's asymptotic evaluation yields  
\begin{equation} \label{eqn:clm2}
\begin{split}
\hgf(\ba + m \bp; z) 
&\approx m^{a+b-c} \cdot \sqrt{2 \pi} \, 
\frac{u(x_0)}{\sqrt{\phi''(x_0)}} \, \Phi(x_0)^m \, m^{-\frac{1}{2}}  
\\[1mm]
&= \frac{2 \sqrt{\pi} \, (2 \sqrt{1-z})^{\, c-a-b-\frac{1}{2} } }{ 
(2m)^{c-a-b+\frac{1}{2}} \, (1+\sqrt{1-z})^{2 m + c-1}}.  
\end{split} \tag{$\ref{eqn:clm}'$}
\end{equation}
For $n = 2 m$ even, since $y_1^{(0)}(n) = \hgf(\ba + m \bp; z)$, 
formula \eqref{eqn:clm} is a direct consequence of \eqref{eqn:clm2}. 
For $n = 2m+1$ odd, since $y_1^{(0)}(n) = \hgf(\ba + m \bp + \bk; z)$, 
formula \eqref{eqn:clm} is obtained from \eqref{eqn:clm2} by 
replacing $\ba$ with $\ba + \bk$. 
Hence the proposition is proved. \hfill $\Box$ 
\section{Discrete Laplace Method} \label{sec:dlm}
In \cite[\S 5]{EI2} we developed a discrete analogue of Laplace's method 
for a class of hypergeometric sums with a large parameter $n$.  
The assumptions imposed there were unnecessarily too restrictive.  
We are able to relax them to some extent without essential changes 
in the proofs so that the improved results should have broader applicability.  
Consider a sum of the form    
$$
g(n) = \sum_{k= \lceil r_0 n \rceil}^{\lceil r_1 n \rceil -1 } G(k; n) \, z^k, 
\qquad 
G(k; n) := \frac{\prod_{i \in I} 
\vG(\sigma_i k + \lambda_i n + \alpha_i)}{\prod_{j \in J} 
\vG(\tau_j k + \mu_j n + \beta_j) },    
$$  
with an independent variable $z$, where $0 \le r_0 < r_1 \le + \infty$; 
$\sigma_i$, $\tau_j \in \bR^{\times}$; $\lambda_i$,  $\mu_j \in \bR$; 
$\alpha_i$, $\beta_j \in \bC$, with $I$, $J$ being finite sets of indices. 
The cardinality of $I$ is denoted by $|I|$. 
Put 
\begin{equation} \label{eqn:rng}
\rho := z \, 
\frac{\prod_{i \in I} |\sigma_i|^{\sigma_i}}{\prod_{j \in J} |\tau_j|^{\tau_j}}, 
\quad 
\nu := \sum_{i \in I} \lambda_i - \sum_{j \in J} \mu_j, 
\quad  
\gamma := \sum_{i \in I} \alpha_i - \sum_{j \in J} \beta_j + 
\frac{|J|-|I|}{2}. 
\end{equation}
\begin{assumption} \label{ass:dlm} Suppose that $z > 0$ and the following 
four conditions are satisfied.   
\begin{enumerate}
\item Balancedness: $\bsigma = (\sigma_i)$ and $\btau = (\tau_j)$ are 
balanced to the effect that 
$$
\sum_{i \in I} \sigma_i = \sum_{j \in J} \tau_j. 
$$
\item Positivity: all gamma factors in $G(k; n)$ are positive 
to the effect that 
$$
l_i(x) := \sigma_i x + \lambda_i  > 0, \qquad  
m_j(x) := \tau_j x + \mu_j > 0, \qquad 
r_0 < {}^{\forall}x < r_1. 
$$
\item Genericness of parameters: $\bal := (\alpha_i) \times (\beta_j)$ is 
generic to the effect that 
$$
\delta_{*}(n; \bal) := \min \left\{1, \, \prod_{i \in I_{*}} 
\dist \left(\alpha_i^{(*)}(n), \,\, \bZ_{\le 0} + |\sigma_i| \bZ_{\le - *} \right)  
\right\} > 0, \qquad * = 0, 1, 
$$
where $\dist(z, \, Z)$ stands for the distance of a point $z \in \bC$ 
from a set $Z \subset \bC$ and 
$$
I_{*} := \{ i \in I \,:\, l_i(r_{*}) = 0 \}, 
\qquad 
\alpha_i^{(*)}(n) := \alpha_i + \sigma_i(\lceil r_{*} n \rceil - r_{*} n). 
$$
\item Convergence: when $r_1 = +\infty$, the infinite series $g(n)$ 
is absolutely convergent for every $n \gg 1$, which is the case if and 
only if one of the following conditions is satisfied: 
\begin{equation} \label{eqn:conv}
\mathrm{(i)} \quad 0 < \rho < 1; \quad 
\mathrm{(ii)} \quad \rho = 1, \,\, \nu < 0; \quad 
\mathrm{(iii)} \quad \rho = 1, \,\, \nu = 0, \,\, \rRe \, \gamma < -1.  
\end{equation}
\end{enumerate}
\end{assumption}
\begin{remark} \label{rem:ass} 
Three remarks are in order about Assumption \ref{ass:dlm}. 
\begin{enumerate} 
\item Balancedness of $\blambda =(\lambda_i)$ and $\bmu = (\mu_j)$, 
that is, the nullity of $\nu$ was assumed in \cite[\S 5]{EI2}, but this condition 
is not essential and hence removed in this article.  
Another improvement is to allow the existence of an independent variable 
$z$, which was fixed to be one in \cite[\S 5]{EI2}. 
\item If $r_{*}$ is an integer then $\alpha_i^{(*)}(n) = \alpha_i$ and 
so $\delta_{*}(n; \bal)$ is independent of $n$, in which case 
$\delta_{*}(n; \bal)$ is simply denoted by $\delta_{*}(\bal)$. 
This will often be the case in practical applications.   
\item If $r_1 = +\infty$ then the positivity (2) forces $\sigma_i$ and  
$\tau_j$ to be positive. 
Stirling's formula gives 
$$
G(k; n) \, z^k \approx \mathrm{const.} \, k^{\nu n + \gamma} \, \rho^k 
\quad \mbox{as} \quad k \to +\infty,  
$$
where $\mathrm{const}.$ is independent of $k$ (but may depend on $n$). 
This asymptotics readily leads to the convergence conditions (i), (ii), (iii) 
in item (4) of Assumption \ref{ass:dlm}.   
\end{enumerate}
\end{remark}
\par
The multiplicative phase function $\Phi(x)$ and the amplitude function  
$u(x)$ are defined by 
$$
\Phi(x) := z^x \, \frac{\prod_{i \in I} l_i(x)^{l_i(x)}}{\prod_{j \in J} m_j(x)^{m_j(x)}}, 
\qquad 
u(x) := (2 \pi)^{\frac{|I|-|J|}{2}} \frac{ \prod_{i \in I} l_i(x)^{\alpha_i-\frac{1}{2}} }{ 
\prod_{j \in J} m_j(x)^{\beta_j-\frac{1}{2}}}, 
\qquad r_0 < x < r_1.  
$$
If $r_1$ is finite then $\Phi(x)$ extends to a positive continuous function 
on the bounded closed interval $[r_0, \, r_1]$.  
If $r_1 = + \infty$ then the function $\Phi(x)$ admits an asymptotic behavior  
$$
\Phi(x) = c \cdot \rho^x \left\{1 + O\left(x^{-1}\right)\right\} \quad 
\mbox{as} \quad x \to + \infty, \qquad 
c := \frac{\prod_{i \in I} |\sigma_i|^{\lambda_i}}{\prod_{j \in J} |\tau_j|^{\mu_j}} 
> 0,   
$$
so one can put $\Phi(+\infty) := 0$ if $0 < \rho < 1$ and 
$\Phi(+\infty) := c > 0$ if $\rho = 1$. 
Thus under convergence condition (4) in Assumption \ref{ass:dlm} $\Phi(x)$ 
extends to a continuous function on $[r_0, \, + \infty]$, which is positive 
on $[r_0, \, +\infty)$. 
In either case $\Phi(x)$ attains a maximum value on $[r_0, \, r_1]$.  
Let
$$
\Phi_{\max} := \max_{r_0 \le x \le r_1} \Phi(x) > 0, \qquad 
\Mm := \{ x \in [r_0, \, r_1] : \Phi(x) = \Phi_{\max} \}. 
$$   
The additive phase function $\phi(x)$ is defined by $\Phi(x) = \re^{- \phi(x)}$.  
A little calculation shows  
$$
\phi'(x) = \log \frac{\prod_{j \in J} m_j(x)^{\tau_j}}{z \prod_{i \in I} l_i(x)^{\sigma_i}}, 
\qquad
\phi''(x) = \sum_{j \in J} \frac{\tau_j^2}{m_j(x)} - 
\sum_{i \in I} \frac{\sigma_i^2}{l_i(x)}. 
$$
Note that any $x_0 \in \Mm \cap (r_0, \, r_1)$ is a solution to the equation 
$\phi'(x) = 0$ or equivalently, 
$$
\prod_{j \in J} m_j(x)^{\tau_j} - z \prod_{i \in I} l_i(x)^{\sigma_i} = 0. 
$$   
\par
We are able to generalize \cite[Theorem 5.2 and Proposition 5.14]{EI2} 
in the following manner. 
\begin{theorem} \label{thm:dlm} 
Suppose that $\Phi(r_*) < \Phi_{\max}$ for $* = 0$, $1$ and that each 
maximum point $x_0 \in \Mm$ is non-degenerate to the effect that 
$\phi''(x_0) > 0$. 
Then $g(n)$ can be expressed as  
\begin{equation} \label{eqn:dlm1}
g(n) = n^{\gamma + \frac{1}{2}} \left(\frac{n}{\re} \right)^{\nu n} 
\Phi_{\max}^{\, n} \{ C + \Omega(n) \}, 
\qquad 
C := \sqrt{2 \pi} \sum_{x_0 \in \Mm} \frac{u(x_0)}{\sqrt{\phi''(x_0)}},  
\end{equation}
and there exist constants $K > 0$, $\lambda > 1$ and $N \in \bN$ 
such that the error term $\Omega(n)$ satisfies 
\begin{equation} \label{eqn:dlm2}
|\Omega(n)| \le K \left\{ n^{-\frac{1}{2}} + \lambda^{-n} 
\left( \delta_0(n; \bal)^{-1} + \delta_1(n; \bal)^{-1} \right) \right\}, 
\qquad {}^{\forall} n \ge N.  
\end{equation}
\end{theorem}
\begin{proposition} \label{prop:dlm} 
For any $\Psi > \Phi_{\max}$ there exist $K > 0$ and $N \in \bN$ such that 
\begin{equation} \label{eqn:dlm3}
|g(n)| \le K \, (n/\re)^{\nu n} \, \Psi^n  \{ \delta_0(n; \bal)^{-1} 
+ \delta_1(n; \bal)^{-1} \} , \qquad {}^{\forall} n \ge N.   
\end{equation}
\end{proposition}
\begin{remark} \label{rem:dlm} 
The constants $K$ and $N$ in \eqref{eqn:dlm2} and \eqref{eqn:dlm3} can 
be taken uniformly with respect to the parameters $\bal$ in any bounded 
subset of $\bC^I \times \bC^J$ (satisfying $\rRe \, \gamma \le -1-\ve$ 
with a fixed $\ve > 0$ if $r_1 = + \infty$ and case (iii) occurs in 
\eqref{eqn:conv}). 
This remark continues to (1) of Remark \ref{rem:dom}.   
\end{remark}
\par
Note that if $r_1 = +\infty$ then $I_1 = \emptyset$ and hence   
$\delta_1(\bal) = 1$ in estimates \eqref{eqn:dlm2} and \eqref{eqn:dlm3}.     
What is new in Theorem \ref{thm:dlm} and Proposition \ref{prop:dlm} is 
the occurrence of the factor $(n/ \re)^n$ in formulas \eqref{eqn:dlm1} 
and \eqref{eqn:dlm3}. 
The proofs of them are practically the same as those of 
\cite[Theorem 5.2 and Proposition 5.14]{EI2}. 
The only difference lies in the manipulation of the function  
$$
H(x; n) \, z^k := \frac{\prod_{i \in I} \vG( l_i(x) n + \alpha_i)}{\prod_{j \in J} 
\vG( m_j(x) n + \beta_j) } \, z^k,  \qquad r_0 < x < r_1.  
$$  
Indeed an application of Stirling's formula to $H(x; n)$ shows that 
as $n \to \infty$, 
$$
H(x; n) \, z^k \approx n^{\gamma} \, u(x) \, \Phi(x)^n \, (n/\re)^{ 
\{ \sum_{i \in I} l_i(x) - \sum_{j \in J} m_j(x) \} n} = 
n^{\gamma} \, u(x) \, \Phi(x)^n \, (n/\re)^{\nu n}, 
$$
by the balancedness (1) in Assumption \ref{ass:dlm} and the 
definition of $\nu$.  
See the proof of \cite[Lemma 5.3]{EI2}, where $\blambda$ and 
$\bmu$ are balanced, i.e., $\nu = 0$, so the factor 
$(n/\re)^{\nu n}$ does not occur.   
\section{Some Examples} \label{sec:ex}
We illustrate Theorem \ref{thm:dlm} and Proposition \ref{prop:dlm}  
by a couple of examples.  
They will be applied to asymptotic analysis of the truncation error for  
Gauss's continued fraction in \S \ref{sec:dom} and \S \ref{sec:error}. 
In this section $\rho$, $\nu$ and $\gamma$ are the ones defined in 
\eqref{eqn:rng} and other notations in \S \ref{sec:dlm} are also retained.     
\begin{example} \label{ex:1}
For $a$, $b$, $c \in \bC$ and $z > 0$ we consider the infinite sum 
\begin{equation} \label{eqn:ex1}
g_1(n; a, b, c; z) := \sum_{k=0}^{\infty} 
\frac{\vG(k+n+a) \, \vG(k+n+b)}{\vG(k+1) \, \vG(k+c)} \, z^k.  
\end{equation}
Note that $r_0 = 0$, $r_1 = +\infty$, $\nu = 2$, $\gamma = a+b-c-1$, 
$l_1(x) = l_2(x) = x+1$, $m_1(x) = m_2(x) = x$,   
\begin{alignat*}{3}
\Phi(x) &= \frac{ z^x \, (x+1)^{2(x+1)}}{ x^{2x}}, 
\qquad &  
u(x) &= \frac{(x+1)^{a+b-1}}{x^{c}},   
\qquad &   & 
\\[1mm]  
\phi'(x) &= 2 \log \frac{x}{\sqrt{z} (x+1)}, 
\qquad & 
\phi''(x) &= \frac{2}{x(x+1)} > 0, 
\qquad & x &\in (0, \, +\infty).  
\end{alignat*}
Since $\rho = z$, the convergence condition is just $0 < z < 1$. 
Under this condition the equation $\phi'(x) = 0$ has a unique solution 
$x_0 = \sqrt{z}/(1-\sqrt{z})$ in $(0, \, +\infty)$.   
Observe that 
\begin{alignat*}{2}
\Phi(x_0) &= (1 - \sqrt{z})^{-2}, \qquad &  
u(x_0) &= z^{-\frac{c}{2}} (1-\sqrt{z})^{c-a-b+1}, \\[1mm]
\phi''(x_0) &= 2 z^{-\frac{1}{2}} (1-\sqrt{z})^2, \qquad & 
\delta_0(\bal) &= \delta_1(\bal) = 1. 
\end{alignat*}
Hence for $z \in (0, \, 1)$ Theorem \ref{thm:dlm} leads to an 
asymptotic representation  
\begin{equation} \label{eqn:ex1a}
g_1(n; a, b, c; z) \sim \frac{ \sqrt{\pi} }{\sqrt{z}^{\, c-\frac{1}{2}}}  
\left(\frac{n}{\re} \right)^{2 n}  
\frac{(1 - \sqrt{z})^{c-a-b-2n}}{ n^{c-a-b+\frac{1}{2}}}.  
\end{equation}
\end{example}
\begin{example} \label{ex:2} 
For $a$, $b$, $c \in \bC$ and $z > 0$ we consider the sum 
$$
g_2(n; a, b, c; z) := \sum_{k=0}^{n-1} 
\frac{\vG(k+n+a)}{ \vG(k+1) \, \vG(k+c) \, \vG(-k+n+b) } \, z^k.  
$$
Note that $r_0 = 0$, $r_1 = 1$, $\nu = 0$, $\gamma = a-b-c$, 
$l_1(x) = x+1$, $m_1(x) = m_2(x) = x$, $m_3(x) = 1-x$,  
\begin{alignat*}{3}
\Phi(x) &= \frac{ z^x \, (x+1)^{x+1}}{ x^{2x} (1-x)^{1-x}}, 
\quad & 
u(x) &= \frac{1}{2 \pi} \cdot \frac{(x+1)^{a-\frac{1}{2}}}{x^c (1-x)^{b-\frac{1}{2}}},  
\qquad & & \\[1mm]  
\phi'(x) &= \log \frac{x^2}{z (1-x^2)}, 
\qquad & 
\phi''(x) &= \frac{2}{x(1-x^2)} > 0, 
\quad & x &\in (0, \, 1).  
\end{alignat*}
The equation $\phi'(x) = 0$ has a unique solution 
$x_0 = \sqrt{z (1+z)^{-1}}$ in $(0, \, 1)$.  
Observe that 
\begin{alignat*}{2}
\Phi(x_0) &= (\sqrt{z} + \sqrt{z+1})^2, \qquad &  
u(x_0) &= (2 \pi)^{-1} z^{-\frac{c}{2}} (z+1)^{\frac{b+c-a}{2}} 
(\sqrt{z} + \sqrt{z+1})^{a+b-1}, \\[1mm]
\phi''(x_0) &= 2 z^{-\frac{1}{2}} (z + 1)^{\frac{3}{2}}, \qquad & 
\delta_0(\bal) &= \delta_1(\bal) = 1. 
\end{alignat*}
Hence for $z > 0$ Theorem \ref{thm:dlm} leads to an asymptotic 
representation 
\begin{equation} \label{eqn:ex2a}
g_2(n; a, b, c; z) \sim 
\frac{\sqrt{z+1}^{\, b+c-a - \frac{3}{2}}}{2 \sqrt{\pi} \sqrt{z}^{\, c-\frac{1}{2}}} 
\cdot 
\frac{  (\sqrt{z} + \sqrt{z+1})^{2 n+a+b-1} }{ n^{b+c-a-\frac{1}{2}}  }. 
\end{equation}
\end{example}
\begin{example} \label{ex:3} 
For $a$, $b$, $c$, $d \in \bC$ and $z > 0$ we consider the infinite series  
$$
g_3(n; a, b, c, d; z) := \sum_{k=n}^{\infty} 
\frac{\vG(2 k+ 2 n+a) \, \vG(2 k-2 n+b)}{\vG(2 k+c) \, \vG(2 k+d)} \, z^k.  
$$
Note that $r_0 = 1$, $r_1 =+\infty$, $\rho = z$, $\nu = 0$ and  
$\gamma = a+b-c-d$, so the convergence condition is either 
$0 < z < 1$ or $z = 1$, $\rRe \, \gamma < -1$, which is assumed from now on.  
Since $l_1(x) = 2(x+1)$, $l_2(x) = 2(x-1)$, $m_1(x) = m_2(x) = 2 x$, we have   
\begin{alignat*}{2}
\Phi(x) &= \frac{ z^x \, (x+1)^{2(x+1)}(x-1)^{2(x-1)}}{ x^{4 x}}, 
\quad & 
u(x) &= 2^{a+b-c-d} \, \frac{(x+1)^{a-\frac{1}{2}} (x-1)^{b-\frac{1}{2}}}{x^{c+d-1} }, 
\\[1mm]
\phi'(x) &= \log \frac{x^4}{z (x^2-1)^2} > 0, \quad & 
\phi''(x) &= -\frac{4}{x(x^2-1)} < 0, \qquad x \in (1, \, +\infty).  
\end{alignat*}
Thus $\Phi(x)$ is strictly decreasing in  
$[1, \, +\infty)$ with maximum $\Phi(1) = 16 z$.  
Observe that 
$$
\delta_0(\bal) = 
\min \{ 1, \, \dist(b, \, \bZ_{\le 0}) \}, 
\qquad \delta_1(\bal) = 1.  
$$
Hence for $0 < z < 1$ or $z = 1$, $\rRe \, \gamma <-1$ 
Proposition \ref{prop:dlm} implies that for any $\Psi > 16 z$ there exist a 
constant $K > 0$ and an integer $N \in \bN$ such that
\begin{equation} \label{eqn:ex3a}
|g_3(n; a, b, c, d; z)| < 
\frac{K \cdot \Psi^n}{\min \{ 1, \, \dist(b, \, \bZ_{\le 0}) \}},  
\qquad {}^{\forall} n \ge N. 
\end{equation}
\end{example}
\section{Decomposition and Sign Changes} \label{sec:pi-sc}
The positivity condition (2) in Assumption \ref{ass:dlm} is not always 
satisfied by a general hypergeometric series. 
To cope with this situation we have to 
discuss how to recover the condition.    
\par
Consider an infinite series of the form 
\begin{equation} \label{eqn:g(n)}
g(n) = \sum_{k=0}^{\infty} G(k; n) \, z^k, \qquad 
G(k; n) := \frac{\prod_{i \in I} \vG(\sigma_i k + n \lambda_i + 
\alpha_i)}{\prod_{j \in J} \vG(\tau_j k + n \mu_j + \beta_j)},  
\end{equation}
where $\sigma_i$, $\tau_j$, $\lambda_i$, $\mu_j \in \bZ$ with 
$\sigma_i \neq 0$ and $\tau_j \ne 0$ for $i \in I$ and $j \in J$. 
Let $r_1 < r_2< \cdots < r_m$ be the distinct positive roots of the product  
$\prod_{i \in I} l_i(x) \prod_{j \in J} m_j(x)$ and put $r_0 := 0$ 
and $r_{m+1} := +\infty$ by convention.  
We decompose the series $g(n)$ into $m+1$ components  
\begin{equation} \label{eqn:gs(n)}
g(n) = \sum_{s=0}^m g_s(n), \qquad 
g_s(n) := \sum_{k= \lceil r_s n \rceil}^{\lceil r_{s+1} n \rceil-1} G(k; n) \, z^k.  
\end{equation}
Since each of the linear functions $l_i(x)$ and $m_j(x)$ is either positive 
everywhere or negative everywhere on each interval 
$\vD_s := (r_s, \, r_{s+1})$, one can define index subsets  
$$
I_s^{\pm} := \{ i \in I \,:\, \mbox{ $l_i(x) \gtrless 0$ on $\vD_s$} \}, 
\qquad 
J_s^{\pm} := \{ j \in J \,:\, \mbox{ $m_j(x) \gtrless 0$ on $\vD_s$} \}.  
$$
The corresponding gamma factors in $G(k; n)$ are said to be {\sl positive} 
or {\sl negative} on $\vD_s$.  
\par
Applying Euler's reflection formula $\vG(x) \, \vG(1-x) = \pi/\sin \pi x$ 
to each negative gamma factor of $G(k;n)$ and taking the assumption 
$\sigma_i$, $\tau_j$, $\lambda_i$, $\mu_j \in \bZ$ into account, we have     
\begin{equation} \label{eqn:gs(n)p}
g_s(n) = \pi^{|I_s^-|-|J_s^-|} \, 
\frac{\prod_{j \in J_s^-} \sin \pi \beta_j}{\prod_{i \in I_s^-} \sin \pi \alpha_i} 
\cdot (-1)^{\nu_s^- n} 
\sum_{k= \lceil r_s n \rceil}^{\lceil r_{s+1} n \rceil-1} G_s(k; n) \, z_s^k, 
\quad z_s := (-1)^{\theta_s^-} z,    
\end{equation}
where $\nu_s^- := \sum_{i \in I_s^-} \lambda_i - \sum_{j \in J_s^-} \mu_j 
\in \bZ$, $\theta_s^- := \sum_{i \in I_s^-} \sigma_i  - 
\sum_{j \in J_s^-} \tau_j \in \bZ$, and  
$$
G_s(k; n) := 
\frac{\prod_{i \in I_s^+} \vG(\sigma_i k + \lambda_i n + \alpha_i) 
\prod_{j \in J_s^-} \vG(-\tau_j k - \mu_j n + 1-\beta_j) }{ 
\prod_{j \in J_s^+} \vG(\tau_j k + \mu_j n + \beta_j)  
\prod_{i \in I_s^-}\vG(- \sigma_i k - \lambda_i n + 1-\alpha_i)}.      
$$  
Notice that all gamma factors in $G_s(k; n)$ are positive on $\vD_s$, 
as desired. 
Proceeding from \eqref{eqn:g(n)} to \eqref{eqn:gs(n)p} via 
\eqref{eqn:gs(n)} is referred to as the procedure of 
{\sl decomposition and sign changes}. 
\par
Let $\kappa_s := |I_s^+|+|J_s^-|-|I_s^-|-|J_s^+|$. 
If $z_s$ is positive then the multiplicative phase function $\Phi_s(x)$ and 
the amplitude function $u_s(x)$ for the sum in \eqref{eqn:gs(n)p} have  representations    
$$
\Phi_s(x) = z_s^x \, 
\frac{\prod_{i \in I} |l_i(x)|^{l_i(x)}}{\prod_{j \in J} |m_j(x)|^{m_j(x)}}, 
\qquad 
u_s(x) = (2 \pi)^{\frac{\kappa_s}{2}} \, 
\frac{\prod_{i \in I} |l_i(x)|^{\alpha_i-\frac{1}{2}}}{\prod_{j \in J} 
|m_j(x)|^{\beta_j-\frac{1}{2}}}, 
\qquad 
x \in \vD_s,  
$$ 
which are independent of $s$ up the the first factors on the right-hand 
sides.    
When $z_s$ is negative, we should make a sign change by dividing the 
sum in \eqref{eqn:gs(n)p} into its {\sl even} and {\sl odd} components, 
where the former is the sum over even $k$'s while the latter is the sum 
over odd $k$'s, so that $z_s^2 = z^2$ becomes a new independent 
variable that is positive. 
This procedure is called the {\sl even-odd decomposition}.  
Here is an example illustrating these procedures. 
\begin{example} \label{ex:dsc} 
For $a$, $b$, $c \in \bC$ and $z < 0$ we consider the infinite sum 
$$
g_4(n; a, b, c; z) := \sum_{k=0}^{\infty} 
\frac{\vG(k+n+a) \, \vG(k-n+b)}{\vG(k+1) \, \vG(k+c)} \, z^k.  
$$
It is absolutely convergent if and only if either $-1< z < 0$ or 
$z = -1$, $\rRe(c-a-b) > 0$, which is assumed from now on.  
The sum decomposes into two components corresponding to $0 \le k \le n-1$ 
and $n \le k < \infty$.  
After the procedure of decomposition and sign changes we have   
$$
g_4(n; a, b, c; z) = \frac{\pi (-1)^n}{\sin \pi b} \cdot  
g_2(n; a, 1-b, c; - z) + h(n; a, b, c; z), 
$$
where $g_2(n; a, b, c; z)$ is defined in Examples \ref{ex:2}, while 
$$
h(n; a, b, c; z) := \sum_{k=2 n}^{\infty} 
\frac{\vG(k+n+a) \, \vG(k-n+b)}{\vG(k+1) \, \vG(k+c)} \, z^k.  
$$
\par
The result \eqref{eqn:ex2a} in Example \ref{ex:2} shows that 
$$
g_2(n; a, 1-b, c; -z) \sim 
\frac{(2 \sqrt{1-z})^{c-a-b-\frac{1}{2}} }{\sqrt{\pi} \, \sqrt{-z}^{\, c-\frac{1}{2}} } 
\cdot  
\frac{(\sqrt{-z} + \sqrt{1-z})^{2 n+a-b} }{ (2n)^{c-a-b+\frac{1}{2} } }. 
$$
According to whether $n = 2m$ or $n = 2m+1$ the even-odd decomposition 
of $h(n; a, b, c; z)$ reads    
\begin{align*}
h(2m; a, b, c; z) &= g_3(m; a, b, 1, c; z^2) + z \, g_3(m; a+1, b+1, 2, c+1; z^2), 
 \\[1mm]
h(2m+1; a, b, c; z) &= 
z \, g_3(m; a+2, b, 2, c+1; z^2) + z^2 \, g_3(m; a+3, b+1, 3, c+2; z^2), 
\end{align*}
where $g_3(n; a, b, c, d; z)$ is defined in Example \ref{ex:3}. 
So the result  \eqref{eqn:ex3a} in this example shows that for any 
$\Psi > 4(-z)$ there exists a constant $K > 0$ and an integer 
$N > 0$ such that 
$$
|h(n; a, b, c; z)| \le \frac{K \cdot \Psi^n}{\min\{1, \, \dist(b, \, \bZ_{\le 0}) \}}, 
\qquad n \ge N, 
$$
whether $n$ is even or odd.  
Taking $\Psi$ so that $4(-z) < \Psi < (\sqrt{-z} + \sqrt{1-z})^2$, we have 
\begin{equation} \label{eqn:g4a}
g_4(n; a, b, c; z) \sim 
\frac{\sqrt{\pi} \, (-1)^n }{ \sin \pi b } \cdot 
\frac{(2 \sqrt{1-z})^{c-a-b-\frac{1}{2}} }{ \sqrt{-z}^{\, c-\frac{1}{2} } } \cdot  
\frac{(\sqrt{-z} + \sqrt{1-z})^{2 n+a-b} }{(2n)^{c-a-b+\frac{1}{2} } }. 
\end{equation}
\end{example}
\section{Dominant Solutions} \label{sec:dom}
According to whether $z \in (0, \, 1)$ or $z \in (-\infty, \, 0)$, 
we take different kinds of dominant solutions to the recurrence equation 
\eqref{eqn:rr}, that is, the solution $y_1^{(1)}(n)$ in the former case and 
a Pfaff transformation of $y_1^{(\infty)}(n)$ in the latter case 
respectively; see Remark \ref{rem:Pfaff} for Pfaff's transformations.       
\begin{lemma} \label{lem:y11} 
For any $z \in (0, \, 1)$ we have  
\begin{equation} \label{eqn:y11} 
y_1^{(1)}(n) 
\sim \frac{\sqrt{\pi} \sin \pi c}{\sin \pi(c-a) \sin \pi(c-b)} 
\cdot \frac{(2 \sqrt{1-z})^{c-a-b-\frac{1}{2}}}{ n^{c-a-b+\frac{1}{2}}  }  
\left(\frac{1+\sqrt{1-z}}{z} \right)^{n+c-1}.  
\end{equation}
\end{lemma}
{\it Proof}. 
From definitions \eqref{eqn:r-frob1} and \eqref{eqn:ex1} we have  
$$
y_1^{(1)}(\ba+m \bp; z) = \chi(\ba+ m \bp) \, 
g_1(m; a, b, a+b-c+1; 1-z),   
$$
where definition \eqref{eqn:m1(a)} and Stirling's formula yields   
\begin{equation} \label{eqn:m1a}
\chi(\ba+ m \bp) \approx \frac{\sin \pi c}{2 \sin \pi(c-a) \sin \pi(c-b)} 
\cdot m^{a+b-2c+1} \left( \frac{m}{\re} \right)^{- 2m}.  
\end{equation}
This together with formula \eqref{eqn:ex1a} in Example \ref{ex:1} leads to  
\begin{equation} \label{eqn:y11d}
y_1^{(1)}(\ba + m \bp) 
\sim \frac{\sqrt{\pi} \sin \pi c}{\sin \pi(c-a) \sin \pi(c-b)} 
\cdot \frac{(2 \sqrt{1-z})^{c-a-b-\frac{1}{2}}}{ (2m)^{c-a-b+\frac{1}{2}}  } 
\left(\frac{1+\sqrt{1-z}}{z} \right)^{2m+c-1}.  
\tag{$\ref{eqn:y11}'$}
\end{equation}
When $n = 2 m$ is even, since $y_1^{(1)}(n) = y_1^{(1)}(\ba + m \bp; z)$, 
formula \eqref{eqn:y11} directly follows from \eqref{eqn:y11d}.  
When $n = 2m+1$ is odd, in view of $y_1^{(1)}(n) = 
y_1^{(1)}(\ba + m \bp + \bk; z)$, formula \eqref{eqn:y11} is obtained from 
\eqref{eqn:y11d} by replacing $\ba$ with $\ba + \bk$. 
Thus the lemma is proved. \hfill $\Box$
\begin{lemma} \label{lem:y1i} 
For any $z \in (-\infty, \, 0)$ we have 
\begin{equation} \label{eqn:y1i} 
y_1^{(\infty)}(n) 
\sim \frac{\sqrt{\pi} \, (-1)^n}{\sin \pi(c-a) \sin \pi(c-b)} 
\cdot \frac{(2 \sqrt{1-z})^{c-a-b-\frac{1}{2}}}{ n^{c-a-b+\frac{1}{2}}  }  
\left(\frac{1+\sqrt{1-z}}{-z} \right)^{n+c-1}.  
\end{equation}
\end{lemma}
{\it Proof}. Recall that $y_1^{(\infty)}(\ba; z)$ is defined by 
\eqref{eqn:r-frob-inf} with \eqref{E9}. 
A Pfaff transformation of it reads  
\begin{equation} \label{eqn:y1ip}
y_1^{(\infty)}(\ba; z) = \frac{\chi(\ba)}{\sin \pi c} \, (1-z)^{-a} 
\hgf(a, c-b; a-b+1; (1-z)^{-1}),  
\end{equation}
which corresponds to formula (11) in \cite[Chap. II, \S 2.8]{Erdelyi}, 
where $\chi(\ba)$ is defined in \eqref{eqn:m1(a)}. 
So 
$$
y_1^{(\infty)}(\ba + m \bp; z) = \frac{\chi(\ba+ m \bp)}{\sin \pi c} \, 
(1-z)^{-m-a} g_1(m; a, c-b; a-b+1; (1-z)^{-1}), 
$$
where $g_1(n; a, b, c; z)$ is defined in \eqref{eqn:ex1}. 
If $z \in (-\infty, \, 0)$ then $(1-z)^{-1} \in (0, \, 1)$, so the result 
\eqref{eqn:ex1a} in Example \ref{ex:1} is applicable. 
It follows from formulas \eqref{eqn:ex1a} and \eqref{eqn:m1a} that 
\begin{equation} \label{eqn:y1id} 
y_1^{(\infty)}(\ba +m \bp; z) 
\sim \frac{\sqrt{\pi}}{\sin \pi(c-a) \sin \pi(c-b)} 
\cdot \frac{(2 \sqrt{1-z})^{c-a-b-\frac{1}{2}}}{ (2 m)^{c-a-b+\frac{1}{2}}  }  
\left(\frac{1+\sqrt{1-z}}{-z} \right)^{2 m+c-1}. 
\tag{$\ref{eqn:y1i}'$}  
\end{equation}
When $n = 2 m$ is even, since $y_1^{(\infty)}(n) = 
y_1^{(\infty)}(\ba + m \bp; z)$, 
formula \eqref{eqn:y1i} directly follows from \eqref{eqn:y1id}.  
When $n = 2m+1$ is odd, in view of $y_1^{(\infty)}(n) = 
y_1^{(\infty)}(\ba + m \bp + \bk; z)$, formula \eqref{eqn:y1i} is 
obtained from \eqref{eqn:y1id} by replacing $\ba$ with $\ba + \bk$. 
Thus the lemma is proved. \hfill $\Box$
\begin{remark} \label{rem:dom} 
Two remarks are in order about Lemmas \ref{lem:y11} and \ref{lem:y1i}. 
\begin{enumerate} 
\item Due to Remark \ref{rem:dlm} the relations $\sim$ in \eqref{eqn:y11} 
and \eqref{eqn:y1i} are compatible with the specialization procedure of 
letting $a \to 0$ followed by the substitution $c \mapsto c-1$.  
\item We wonder whether in the proof of Lemma \ref{lem:y11} a Pfaff 
transformation of $y_1^{(1)}(n)$ could be employed instead of itself.   
A Pfaff transformation of $y_1^{(1)}(\ba; z)$ in \eqref{eqn:r-frob1} 
reads 
$$
y_1^{(1)}(\ba; z) = \frac{\sin \pi c}{\sin \pi(c-b)} \cdot 
\frac{\vG(b)}{\vG(c-b)} \, z^{-a} 
\hgf(a, a-c+1; a+b-c+1; 1-z^{-1}),  
$$
which is a rescaled version of formula (7) in 
Erd\'elyi \cite[Chap. II, \S 2.8]{Erdelyi}.  
So we have  
\begin{align*}
y_1^{(1)}(\ba + m \bp; z) 
&= \frac{(-1)^m \sin \pi c}{\sin \pi (c-b)} \cdot 
\frac{\vG(m+b)}{\vG(m+c-b)}  \\
&\phantom{=} \times z^{-m-a} g_4(m; a, a-c+1, a+b-c+1; 1-z^{-1}), 
\end{align*}
where $g_4(n; a, b, c; z)$ is defined in Example \ref{ex:dsc}. 
Note that $\vG(m+b)/\vG(m+c-b) \approx m^{2b-c}$ by Stirling's 
formula.     
Thus the result \eqref{eqn:g4a} in Example \ref{ex:dsc} implies   
formula \eqref{eqn:y11d}, but unfortunately it is valid only for 
$z \in (1/2, \, 1)$ not for all $z \in (0, \, 1)$.  
Similarly, in the proof of Lemma \ref{lem:y1i} the use of 
$y_1^{(\infty)}(n)$ itself in stead of its Pfaff transformation   
leads to formula \eqref{eqn:y1id}, but it is valid only for 
$z \in (-\infty, \, -1)$ not for all $z \in (-\infty, \, 0)$. 
\end{enumerate}     
\end{remark}
\section{Casoratian and Error Estimates} \label{sec:error}
To use error estimate \eqref{eqn:error} we have to evaluate 
the Casoratian $\omega(0)$.  
Let $\bk := (1, 0; 1)$ and 
\begin{align*}
\omega^{(0)}(\ba; z) &:= y_1^{(0)}(\ba; z) \, y_2^{(0)}(\ba+\bk; z) - 
y_1^{(0)}(\ba+\bk; z) \, y_2^{(0)}(\ba; z), \\[1mm]
\omega^{(1)}(\ba; z) &:= y_1^{(0)}(\ba; z) \, y_1^{(1)}(\ba+\bk; z) - 
y_1^{(0)}(\ba+\bk; z) \, y_1^{(1)}(\ba; z), \\[1mm]
\omega^{(\infty)}(\ba; z) &:= y_1^{(0)}(\ba; z) \, y_1^{(\infty)}(\ba+\bk; z) - 
y_1^{(0)}(\ba+\bk; z) \, y_1^{(\infty)}(\ba; z).   
\end{align*}
\begin{lemma} \label{lem:casorati} 
We have $\omega^{(0)}(\ba; z) = - \omega^{(1)}(\ba; z)$ and 
\begin{align*}
\omega^{(1)}(\ba; z) 
&= \frac{\pi \sin \pi c}{\sin \pi(c-a) \cdot \sin \pi(c-b)} 
\cdot \frac{\vG(a) \, \vG(b)}{\vG(c-a) \, \vG(c-b+1)} \cdot 
z^{-c} (1-z)^{c-a-b}, \\[1mm]
\omega^{(\infty)}(\ba; z) 
&= - \frac{\pi}{\sin \pi(c-a) \cdot \sin \pi(c-b)} 
\cdot \frac{\vG(a) \, \vG(b)}{\vG(c-a) \, \vG(c-b+1)} \cdot 
(-z)^{-c} (1-z)^{c-a-b}.  
\end{align*}
\end{lemma}
{\it Proof}.  
It follows from connection formula \eqref{E35} that 
$- \omega^{(1)}(\ba; z) = \omega^{(0)}(\ba; z)$.  
Let  
\begin{align*}
W(\ba; z) := y_1^{(0)}(\ba; z) \, y_2^{(0)}(\ba+\1; z) - 
y_1^{(0)}(\ba + \1; z) \, y_2^{(0)}(\ba; z), 
\qquad \1 := (1, 1; 1). 
\end{align*}
As in the proof of \cite[Lemma 2.1, formula (17c)]{EI1} for ${}_3F_2(1)$,  
we can show  
$$
\frac{d}{d z} y_i^{(0)}(\ba; z) = y_i^{(0)}(\ba+\1; z), \qquad i = 1, 2, 
$$
so that $W(\ba; z)$ is the Wronskian of $y_1^{(0)}(\ba; z)$ and $y_2^{(0)}(\ba; z)$.  
A simple calculation yields  
\begin{equation} \label{eqn:W}
W(\ba; z) = \frac{\pi \sin \pi c}{\sin \pi(c-a) \cdot \sin \pi(c-b)} 
\cdot \frac{\vG(a) \, \vG(b)}{\vG(c-a) \, \vG(c-b)} \cdot 
z^{-c} (1-z)^{c-a-b-1}. 
\end{equation}
There is a simultaneous contiguous relations for the six functions in 
\eqref{eqn:r-frob}, 
$$
y(\ba+ \bk; z) 
= \frac{a}{c-b} \, y(\ba; z) - \frac{1-z}{c-b} \, y(\ba+\1; z).  
$$
Using this relation for $y(\ba; z) = y_i^{(0)}(\ba; z)$, $i = 1$, $2$,  
we have 
$$
- \omega^{(1)}(\ba; z) = \omega^{(0)}(\ba; z) = -\frac{1-z}{c-b} \, W(\ba; z). 
$$
This together with formula \eqref{eqn:W} proves the lemma. 
\hfill $\Box$ \par\medskip
Now we are in a position to establish Theorem \ref{thm:Error} and 
Corollary \ref{cor:Error}.  
\par\medskip\noindent
{\it Proof of Theorem $\ref{thm:Error}$}.  
For $z = 0$ formula \eqref{eqn:Error} is trivial and there is nothing to 
discuss.   
\par
For $z \in (0, \, 1)$ we apply the general estimate \eqref{eqn:error} 
to $f(n) = y_1^{(0)}(n)$ and $g(n) = y_1^{(1)}(n)$. 
Note that $f(0) = f_1^{(0)}(\ba; z)$, $g(0) = \chi(\ba) \, f_1^{(1)}(\ba; z)$ 
and $\omega(0) = \omega^{(1)}(\ba; z)$. 
If we put $w := z(1+\sqrt{1-z})^{-2} > 0$, then it follows from 
Proposition \ref{prop:clm} and Lemma \ref{lem:y11} that  
$$
h(n) = \frac{y_1^{(0)}(n+2)}{y_1^{(1)}(n+2)} \sim  
2 \, \frac{\sin \pi(c-a) \sin \pi(c-b)}{\sin \pi c} \cdot 
w^{n+c+1}. 
$$
Using this formula, various definitions in \S \ref{sec:scr},   
the first formula in Lemma \ref{lem:casorati} as well as the 
recursion formula for the gamma function, we obtain   
\begin{subequations} \label{eqn:y11e}
\begin{align}
\frac{c}{a} \cdot
\frac{\omega(0) \cdot h(n)}{f(0)^2} 
&\sim \frac{2 \pi}{\hgF(\ba; z)^2} \cdot 
\frac{\vG(c) \vG(c+1)}{\vG(a+1) \vG(b) \vG(c-a) \vG(c-b+1)} 
\cdot \frac{z(1-z)^{c-a-b} \, w^n}{(1+\sqrt{1-z})^{2(c+1)}}, 
\label{eqn:y11e1}
\\[2mm] 
\frac{g(0) \cdot h(n)}{f(0)} 
& \sim \frac{\hgF(a, b; a+b-c+1; 1-z)}{\hgF(\ba; z) \, \vG(a+b-c+1)} 
\cdot \frac{2 \pi \, \vG(c) \,  
w^{n+c+1} }{\vG(c-a) \vG(c-b)},    
\label{eqn:y11e2}
\end{align}
\end{subequations} 
where the left-hand side of \eqref{eqn:y11e2} is regular except 
at the poles of $\vG(c)$ and the zeros of $\hgF(\ba; z)$. 
Formula \eqref{eqn:Error} is then derived by combining 
\eqref{eqn:en}, \eqref{eqn:error} and \eqref{eqn:y11e}.  
\par
For $z \in (-\infty, \, 0)$ we apply estimate \eqref{eqn:error} to 
the case where $f(n) = y_1^{(0)}(n)$ and $g(n)$ is a Pfaff transform  
of $y_1^{(\infty)}(n)$.  
Note that $f(0) = f_1^{(0)}(\ba; z)$, $g(0) = y_1^{(\infty)}(\ba; z)$ 
given by \eqref{eqn:y1ip} and $\omega(0) = \omega^{(\infty)}(\ba; z)$. 
Again we put  $w := z(1+\sqrt{1-z})^{-2}$, which is negative this time.     
It follows from Proposition \ref{prop:clm} and Lemma \ref{lem:y1i} that  
$$
h(n) = \frac{y_1^{(0)}(n+2)}{y_1^{(\infty)}(n+2)} \sim  
2 \, \sin \pi(c-a) \sin \pi(c-b) \cdot (-w)^{c+1} w^n.  
$$
Using this formula, various definitions in \S \ref{sec:scr},   
the second formula in Lemma \ref{lem:casorati} as well as the recursion 
formula for the gamma function, we obtain the 
same formula as \eqref{eqn:y11e1} and
\begin{equation} \label{eqn:y1ie}
\frac{g(0) \cdot h(n)}{f(0)} 
\sim  \frac{\hgF(a, c-b; a-b+1; (1-z)^{-1})}{\hgF(\ba; z) \, \vG(a-b+1)} 
\cdot \frac{2 \pi \, \vG(c)}{\vG(b) \vG(c-a)} 
\cdot \frac{(-w)^{c+1} \, w^n}{(1-z)^a }, 
\tag{$\ref{eqn:y11e}$c}
\end{equation}
where the left-hand side of \eqref{eqn:y1ie} is regular except at  
the poles of $\vG(c)$ and the zeros of $\hgF(\ba; z)$. 
Formula \eqref{eqn:Error} is then derived from \eqref{eqn:en}, 
\eqref{eqn:error}, \eqref{eqn:y11e1} and \eqref{eqn:y1ie} as well as 
the reflection formula for the gamma function. 
\hfill $\Box$ \par\medskip\noindent
{\it Proof of Corollary $\ref{cor:Error}$}. 
By item (1) of Remark \ref{rem:dom} the three relations $\sim$ in 
\eqref{eqn:y11e} are compatible with the specialization procedure of 
letting $a \to 0$ followed by the substitution $c \mapsto c-1$. 
Through this procedure the right-hand sides of \eqref{eqn:y11e} 
change in the following manner. 
\begin{align*}
\mbox{RHS of \eqref{eqn:y11e1}} \,  
& \longmapsto \, \frac{2 \pi \vG(c)}{\vG(b) \vG(c-b)} \cdot 
\frac{z(1-z)^{c-b-1}}{(1+\sqrt{1-z})^2} \cdot w^n,  \\[2mm]
\mbox{RHS of \eqref{eqn:y11e2}} \, 
& \longmapsto \,  -2 \sin \pi(c-b) \cdot w^n, \\[1mm] 
\mbox{RHS of \eqref{eqn:y1ie}} \,   
& \longmapsto \, 2 \sin \pi b \cdot (-w)^c w^n. 
\end{align*}
Note that the latter two expressions are regular in $\bb = (b; c)$ and 
hence cause no trouble in applying formulas \eqref{eqn:en} and 
\eqref{eqn:error}. 
Now formula \eqref{eqn:Error2} follows from \eqref{eqn:Error} readily. 
\hfill $\Box$
 

\begin{thebibliography}{99} 
\bibitem{AAR}{G.E.~Andrews, R.~Askey and R.~Roy}, 
{\sl Special Functions}, 
Cambridge Univ. Press, Cambridge, 1999. 
\bibitem{BCP}{J.M.~Borwein, K.-K.S.~Choi and W.~Pigulla}, 
{\sl Continued fractions of tails of hypergeometric series}, 
Amer. Math. Monthly {\bf 112} (2005), no. 6, 493--501. 
\bibitem{CCD}{M.~Colman and A.~Cuyt and J.~Van Deun}, 
{\sl Validated computation of certain hypergeometric functions}, 
ACM Trans. Math. Softw. {\bf 38} (2011), no. 2, Article 11, 20 pages.  
\bibitem{EI1}{A.~Ebisu and K.~Iwasaki}, 
{\sl Three-term relations for ${}_3F_2(1)$}, 
J. Math. Anal. Appl. {\bf 463} (2018), no. 2, 593--610. 
\bibitem{EI2}{A.~Ebisu and K.~Iwasaki}, 
{\sl Contiguous relations, Laplace's methods, and continued fraction for 
${}_3F_2(1)$}, Ramanujan J. (2018) {\tt DOI:10.1007/s11139-018-0039-2}. 
55 pages.  
\bibitem{Erdelyi}{A.~Erd\'elyi et. al.},  
{\sl {H}igher {T}ranscendental {F}unctions, Vol. I}, 
McGraw-Hill, New York, 1953. 
\bibitem{JT}{W.B.~Jones and W.~Thron}, 
{\sl Continued {F}ractions. {A}nalytic {T}heory and {A}pplications}, 
Addison-Wesley, Reading, MA, 1980.   
\end{thebibliography}
\end{document}